\documentclass[a4paper,10pt]{amsart}
\usepackage{amsmath}
\usepackage[T1]{fontenc}
\usepackage{amssymb}
\usepackage{amsthm}
\newcommand{\md}{{\rm d}}

\newcommand{\e}{\varepsilon}
\newtheorem{thm}{Theorem}[section]
\newtheorem{prop}[thm]{Proposition}
\newtheorem{lem}[thm]{Lemma}

\newtheorem{definition}[thm]{Definition}

\newtheorem{rem}[thm]{Remark}
\usepackage{color}

\numberwithin{equation}{section}

\title[G-L energy with
discontinuous constraint] {Magnetic vortices for a Ginzburg-Landau
type energy with discontinuous constraint.~II}

\author[H. Aydi]{Hassen Aydi}
\address{H. Aydi\newline
Universit{\'e} de Monastir\\ Institut sup{\'e}rieur d'informatique de
Mahdia\\ Km 4, R{\'e}jiche,  5121 Mahdia,  Tunisie}
\email{hassen.aydi@isima.rnu.tn}

\author[A. Kachmar]{Ayman Kachmar}
\address{A. Kachmar\newline
Universit{\'e} Paris-Sud\\ D{\'e}partement de math{\'e}matique\\B{\^a}t.
425\\F-91405 Orsay} \email{ayman.kachmar@math.u-psud.fr}

\begin{document}
\maketitle

\begin{abstract}
We study vortex nucleation for minimizers of a Ginzburg-Landau
energy with discontinuous constraint. For applied magnetic fields
comparable with the  first critical field of vortex nucleation, we
determine the limiting vorticities.
\end{abstract}

\section{Introduction and main results}

In the framework of the Ginzburg-Landau theory, it is proposed to
model the energy of an inhomogeneous superconducting sample by
means of the following functional (see \cite{Chetal, Ru})
\begin{equation}\label{V-EGL}
\mathcal G_{\varepsilon,H}(\psi,A)=\int_{\Omega}\left(
|(\nabla-iA)\psi|^2+\frac1{2\varepsilon^2}(p(x)-|\psi|^2)^2+ |{\rm
curl}\,A-H|^2\right)\,\md x.
\end{equation}
Here, $\Omega\subset\mathbb R^2$ is the 2-D cross section of the
superconducting sample, assumed to occupy a cylinder of infinite
height. The complex-valued function $\psi\in H^1(\Omega;\mathbb
C)$ is called the `order parameter', whose modulus $|\psi|^2$
measures the density of the superconducting electron Cooper pairs
(hence $\psi\equiv0$ corresponds to a normal state), and the real
vector field $ A=(A_1,A_2)$ is called the `magnetic potential',
such that the induced magnetic field in the
sample corresponds to ${\rm curl}\, A$.\\
The functional (\ref{V-EGL}) depends on many parameters:
$\frac1\varepsilon=\kappa$ is a characteristic of the
superconducting sample (a temperature independent quantity),
$H\geq0$ is the intensity of the applied magnetic field (assumed
constant and parallel to the axis of the superconducting sample),
$p(x)$ is a positive function modeling the impurities in the
sample, whose values are temperature dependent. The positive sign
of the function $p$ means that the temperature remains below the
critical
temperature of the superconducting sample.\\
In this paper, the function $p$ is a step function. We take $S_1$ an
open smooth set  such that $\overline S_1\subset\Omega$,
$S_2=\Omega\setminus \overline S_1$, and
\begin{equation}\label{p(x)}
p(x)=\left\{\begin{array}{l} 1\quad{\rm if}~x\in S_1\,,\\
a\quad{\rm if}~x\in S_2\,,
\end{array}\right.
\end{equation}
where $a\in\mathbb R_+\setminus\{1\} $ is a given
constant.\\
The above choice of $p$ has  two physical interpretations (see
\cite{Kach-v, Ru}):
\begin{itemize}
\item $S_1$ and $S_2$ correspond to two  superconducting samples
with different critical temperatures; \item The superconducting
sample $\Omega$ is subject to two different temperatures in the
regions $S_1$ and $S_2$, which may happen by cold or heat working
$S_2$.
\end{itemize}
Minimization of  the functional (\ref{V-EGL}) will take place in
the space
$$\mathcal H=H^1(\Omega;\mathbb C)\times
H^1(\Omega;\mathbb R^2)\,.$$ That is we do not assume {\it a priori}
boundary conditions for admissible configurations, but minimizers
satisfy {\it natural boundary conditions}. We study nucleation of
vortices as the applied magnetic field varies in such a manner that,
$$\lim_{\varepsilon\to0}\frac{H}{|\ln\varepsilon|}=\lambda\,,\quad\lambda\in\mathbb R_+\,,$$
and we obtain that their behavior is strongly dependent on the
parameter $a$, see Theorems~\ref{thm1} and \ref{thm2} below.\\
The singularity of the potential $p$ causes a singularity in the
energy as $\varepsilon\to0$ of the order of $\varepsilon^{-1}$,
whereas the energy due to the presence of $n$ vortices typically
diverges like $n\ln\e$. In order to separate the different
singularities of the energy, we let $u_\e$ be the (unique) positive
minimizer of (\ref{V-EGL}) when $H=0$ (see
Theorem~\ref{V-thm-kach3}), the energy of $u_\e$ being of the order
$\e^{-1}$. Then, if $(\psi,A)$ is a minimizing configuration of
(\ref{V-EGL}), it holds that (see Lemma~\ref{V-lem-psi<u}),
$$\mathcal G_{\e,H}(\psi,H)= \mathcal G_{\e,0}(u_\e,0)+\mathcal
F_{\e,H}\left(\frac{\psi}{u_\e},A\right)\,,$$ and the configuration
$\left(\displaystyle\frac{\psi}{u_\e},A\right)$ minimizes the
functional $\mathcal F_{\e,H}$ introduced below,
\begin{equation}\label{reducedfunctional*}
\mathcal
F_{\e,H}(\varphi,A)=\int_{\Omega}\left(u_\e^2|(\nabla-iA)\varphi|^2+\frac{u_\e^4}{2\e^2}(1-|\varphi|^2)^2
+|{\rm curl}\,A-H|^2\right)\,\md x\,.
\end{equation}
The leading order term of the minimizing energy of
(\ref{reducedfunctional*}) will be described by means of an
auxiliary  energy introduced in (\ref{eq-E-lambda}) below. Let
$\mathcal M(\Omega)$ be the space of bounded Radon measures on
$\Omega$, i.e. the topological dual  of $C_0^0(\Omega)$. A measure
$\mu\in\mathcal M(\Omega)$ can be represented canonically as a
difference of two positive measures, $\mu=\mu_+-\mu_-$. The {\it
total variation} and the {\it norm} of $\mu$, denoted respectively
by $|\mu|$ and $\|\mu\|$, are by definition
$|\mu|=\mu_++\mu_-$ and $\|\mu\|=|\mu|(\Omega)$.\\
Given $\lambda>0$, we introduce an energy $E_\lambda$ defined on
$\mathcal M(\Omega)\cap H^{-1}(\Omega)$ as follows. For
$\mu\in\mathcal M(\Omega)\cap H^{-1}(\Omega)$, let $h_\mu\in
H^1(\Omega)$ be the solution of
\begin{equation}\label{eq-h-mu}
\left\{
\begin{array}{l}
-{\rm div}\left(\frac1{p(x)}\nabla h_\mu\right)+h_\mu=\mu\quad{\rm
  in}~
\Omega\,,\\
h_\mu=1\quad{\rm on}~\partial\Omega\,,
\end{array}\right.
\end{equation}
where $p$ is introduced in (\ref{p(x)}). Now $E_\lambda(\mu)$ is
by definition
\begin{equation}\label{eq-E-lambda}
E_\lambda(\mu)=\frac1{\lambda}\int_\Omega p(x)|\mu|\,\md
x+\int_\Omega \left(\frac1{p(x)}|\nabla
h_\mu|^2+|h_\mu-1|^2\right)\,\md x\,.
\end{equation}

\begin{thm}\label{thm1}
Given $\lambda>0$, assume that
$$\lim_{\varepsilon\to0}\frac{H}{|\ln\varepsilon|}=\lambda\,,$$
then
$$\frac{\mathcal F_{\varepsilon,H}}{H^2}\to E_\lambda$$
in the sense of $\Gamma$-convergence. Here the functional $\mathcal
F_{\varepsilon,H}$ is defined in (\ref{reducedfunctional*}) above.
\end{thm}

The convergence in Theorem~\ref{thm1} is precisely described in
Propositions~\ref{prop-upperbound} and \ref{prop-lowerbound} below.\\
Minimizers of (\ref{eq-E-lambda}) can be characterized by means of
minimizers of the following problem,
\begin{equation}\label{eq:dual-en}
\min_{\substack{h\in H^1_0(\Omega)\\-{\rm
div}\left(\frac1{p(x)}\nabla h\right)+h\in\mathcal
M(\Omega)}}\int_\Omega\left(\frac{p(x)}{\lambda}\left|-{\rm
div}\left(\frac1{p(x)}\nabla h\right)+h+1 \right|+\frac{|\nabla
h|^2}{p(x)}+|h|^2\right)\,\md x\,.
\end{equation}
The above functional being strictly convex and continuous, it admits
a unique minimizer, and so the functional $E_\lambda$.\\
Therefore, as a corollary of Theorem~\ref{thm1}, we may describe the
limiting vorticity measure in terms of the minimizer of the limiting
energy $E_\lambda$.

\begin{thm}\label{thm2}
Under the hypothesis of Theorem~\ref{thm1}, if
$(\varphi_\varepsilon, A_\varepsilon)$ is a minimizer of
(\ref{reducedfunctional*}), then, denoting by
$$h_\varepsilon={\rm curl}\,A_\varepsilon\,,\quad
\mu_\varepsilon=h_\varepsilon+
{\rm
  curl}(i\varphi_\varepsilon\,,\,(\nabla-iA_\varepsilon)\varphi
_\varepsilon)\,,$$ the `induced magnetic field' and `vorticity
measure' respectively, the following convergence hold,
\begin{eqnarray}
&&\frac{\mu_\varepsilon}{H}\to \mu_*\quad{\rm in~}
 \left(C^{0,\gamma}(\Omega)\right)^*\quad{\rm
 for~all~}\gamma\in(0,1)\,,\\
&&\frac{h_\varepsilon}{H}\to h_{\mu_*} \quad{\rm
 weakly~in~}H^1_1(\Omega)~{\rm and~strongly~
 in~}W^{1,p}(\Omega)\,,~\forall~p<2\,.
\end{eqnarray}
Here $\mu_*$ is the unique minimizer of $E_\lambda$.
\end{thm}

\subsection*{Sketch of Proof}\ \\
Let us briefly describe the main points of the proofs of
Theorems~\ref{thm1} and \ref{thm2}, and thus explain what stands
behind their
statements.\\
The starting point is the analysis of minimizers of (\ref{V-EGL})
when $H=0$. In this case, (\ref{V-EGL}) has, up to a gauge
transformation, a unique minimizer $(u_\varepsilon,0)$ where
$u_\varepsilon$ is a positive real-valued function. The asymptotic
decay of $u_\varepsilon$ as $\varepsilon\to0$ is obtained in
Lemma~\ref{interior}.\\
When $H>0$, let $(\psi,A)$ be a minimizer of (\ref{V-EGL}). Inspired
by Lassoued-Mironescu \cite{LaMi}, we introduce a {\it normalized
density}\footnote{Notice that $\varphi$
  and $\psi$ have the same vortices.}
$$\varphi=\frac{\psi}{u_\varepsilon}\,.$$
Then $|\varphi|\leq1$ and we are led to the analysis of the
functional $\mathcal F_{\e,H}(\varphi,A)$. The main difficulty stems
from the boundary layer behavior of the weight $u_\e$, which has
rapid oscillations comparable to $\e^{-1}$ near the boundary of
$S_1$. The challenge is then to construct a test configuration with
vortices localized near $\partial S_1$ and having the right amount
of energy. Technically, one does that via a Green's potential $G_\e$
which we were not able to give a good control of it, see
Lemma~\ref{Green}. To overcome this difficulty, we construct
vortices situated at a distance $\frac{\ln|\ln\e|}{|\ln\e|}$ away
from $\partial S_1$. Those vortices being expected to lie on a union
$\gamma$ of  closed curves, a nice $L^1$ bound can be shown to hold
for $G_\e(x,x)$ on $\gamma$. Then, using a result of \cite{Kach-v},
recalled in Lemma~\ref{integration}, this $L^1$ bound provides us
with a family of well separated points serving as the centers of
vortices for the test configurations, see Proposition~\ref{prop-vp}.
We use then suitable methods to estimate from above the energy of
the constructed configuration, yielding in
Proposition~\ref{prop-upperbound} the first part of
Theorem~\ref{thm1}. This rough
analysis will be carried out through Section~\ref{Section:UB}.\\
In Section~\ref{Sec:LB}, we establish a lower bound of the energy,
see Proposition~\ref{prop-lowerbound}. This yields the
second part of Theorem~\ref{thm1}.\\
The proof of Theorem~\ref{thm2} is a by-product of the proofs of
Propositions~\ref{prop-upperbound} and \ref{prop-lowerbound}, see
Remark~\ref{rem:conc}. Finally, we close the paper with a discussion
concerning  the minimization of the energy (\ref{eq-E-lambda}).

\paragraph{\it Remarks on the notation.}
\begin{itemize}
\item
The letters $C,\widetilde C, M,$ etc. will  denote positive
constants independent of $\varepsilon$.
\item For $n\in\mathbb N$
and $X\subset\mathbb R^n$, $|X|$ denotes the Lebesgue measure of
$X$. $B(x,r)$ denotes the open ball in $\mathbb R^n$ of radius $r$
and center $x$.
\item $(\cdot,\cdot)$ denotes the scalar product in
$\mathbb C$ when identified with $\mathbb R^2$.
\item $\mathcal
F_{\e,H}(\varphi,A,U)$ means that the energy density of
$(\varphi,A)$ is integrated only on $U\subset\Omega$.
\item For $\alpha\in(0,1)$ sufficiently small,
$S_1^\alpha=\{x\in S_1~:~{\rm dist}(x,\partial S_1)\geq \alpha\}$
and $S_2^\alpha=\{x\in S_2~:~{\rm dist}(x,\partial
S_1)\geq\alpha\}$.
\item For two positive
functions $a(\varepsilon)$ and $b(\varepsilon)$, we write
$a(\varepsilon)\ll b(\varepsilon)$ as $\varepsilon\to0$ to mean that
$\displaystyle\lim_{\varepsilon\to0}\frac{a(\varepsilon)}{b(\varepsilon)}=0$.
\end{itemize}

\section{Preliminary analysis of minimizers}
\label{V-Sec-EnH=0}

\subsection{The case without applied magnetic field}
This section is devoted to an analysis for minimizers of
(\ref{V-EGL}) when the applied magnetic field $H=0$. We follow
closely similar results obtained in
\cite{kach3}.\\
We keep the notation introduced in Section~1. Upon taking $A=0$
and $H=0$ in (\ref{V-EGL}), one is led to introduce the functional
\begin{equation}\label{V-EnH=0}
\mathcal G_\varepsilon(u):=\int_\Omega\left(|\nabla u|^2
+\frac{1}{2\varepsilon^2}(p(x)-u^2)^2\right)\,\md x\,,
\end{equation}
defined for functions in $H^1(\Omega;\mathbb R)$.\\
We introduce
\begin{equation}\label{V-C0}
C_0(\varepsilon)=\inf_{u\in H^1(\Omega;\mathbb R)} \mathcal
G_\varepsilon(u)\,.
\end{equation}

The next theorem is an analogue of Theorem~1.1 in \cite{kach3}.

\begin{thm}\label{V-thm-kach3}
Given $a\in\mathbb R_+\setminus\{1\}$, there exists
$\varepsilon_0$ such that for all
$\varepsilon\in]0,\varepsilon_0[$, the functional (\ref{V-EnH=0})
admits in $H^1(\Omega;\mathbb R)$ a minimizer $u_\varepsilon\in
C^2(\overline{S_1})\cup C^2(\overline{S_2})$ such that
$$\min(1,\sqrt{a}\,)<u_\varepsilon<\max(1,\sqrt{a}\,)\quad {\rm in}~\overline{\Omega}.$$
If $H=0$, minimizers of (\ref{V-EGL}) are gauge equivalent to the
state $(u_\varepsilon,0)$.
\end{thm}

We state also some estimates, taken from
\cite[Proposition~5.1]{kach3}, that describe the decay of
$u_\varepsilon$ away from the boundary of $S_1$.

\begin{lem}\label{interior}
Let $k\in\mathbb N$. There exist positive constants
$\varepsilon_0$, $\delta$ and $C$ such that, for all
$\varepsilon\in]0,\varepsilon_0]$,
\begin{equation}\label{interiordecay}
\left\|(1-u_\varepsilon)\exp\left(\frac{\delta{\rm dist}(x,\partial
S_1)}\varepsilon\right)\right\|_{H^k(S_1)}+\left\|(\sqrt{a}-u_\varepsilon)\exp\left(\frac{\delta
{\rm dist}(x,\partial
S_1)|}\varepsilon\right)\right\|_{H^k(S_2)}\leq\frac{C}{\varepsilon^k}\,.
\end{equation}
\end{lem}

Finally, we mention without proof that the energy
$C_0(\varepsilon)$ (cf. (\ref{V-C0})) has the order of
$\varepsilon^{-1}$, and we refer to the methods in
\cite[Section~6]{kach3} which permit to obtain the leading order
asymptotic expansion
$$C_0(\varepsilon)=\frac{c_1(a)}\varepsilon+c_2(a)+o(1),\quad(\varepsilon\to0),$$
where $c_1(a)$ and $c_2(a)$ are positive explicit constants.

\subsection{The case with magnetic field}\label{V-sec-minimizers}
This section is devoted to a preliminary analysis of the
minimizers of (\ref{V-EGL}) when $H\not=0$. The main point that we
shall show is how to extract the singular term $C_0(\varepsilon)$
(cf. (\ref{V-C0})) from the energy of a minimizer.

Notice that the existence of minimizers is standard starting from
a minimizing sequence (cf. e.g. \cite{Gi}). A standard choice of
gauge permits one to assume that the magnetic potential satisfies
\begin{equation}\label{V-gauge}
{\rm div}\,A=0\quad {\rm in}~\Omega,\quad n\cdot A=0\quad{\rm
on}~\partial\Omega,
\end{equation}
where $n$ is the outward unit normal vector of
$\partial\Omega$.\\
With this choice of gauge, one is able to prove (since the
boundaries of $\Omega$ and $S_1$ are  smooth) that a minimizer
$(\psi,A)$ is in $ C^1(\overline\Omega;\mathbb C)\times
C^1(\overline\Omega;\mathbb R^2)$. One  has also the following
regularity (cf. \cite[Appendix~A]{kach3}),
$$\psi \in C^2(\overline S_1;\mathbb C)\cup C^2(\overline S_2;\mathbb
C),\quad A\in C^2(\overline S_1;\mathbb R^2)\cup C^2(\overline
S_2;\mathbb R^2).$$

The next lemma is inspired from the work of Lassoued-Mironescu
(cf. \cite{LaMi}).

\begin{lem}\label{V-lem-psi<u}
Let $(\psi,A)$ be a minimizer of (\ref{V-EGL}). Then
$0\leq|\psi|\leq u_\varepsilon$ in $\Omega$, where $u_\varepsilon$
is the positive
minimizer of (\ref{V-EnH=0}).\\
Moreover, putting $\varphi=\frac{\psi}{u_\varepsilon}$, then the
energy functional (\ref{V-EGL}) splits in the form~:
\begin{equation}\label{V-splittingEn}
\mathcal G_{\varepsilon,H}(\psi,A)=C_0(\varepsilon)+ \mathcal
F_{\varepsilon,H}(\varphi,A),\end{equation} where
$C_0(\varepsilon)$ has been introduced in (\ref{V-C0}) and the new
functional $\mathcal F_{\varepsilon,H}$ is defined by~:
\begin{equation}\label{V-reducedfunctional}
\mathcal F_{\varepsilon,H}(\varphi,A)=\int_\Omega
\left(u_\varepsilon^2 |(\nabla-iA)\varphi|^2
+\frac1{2\varepsilon^2} u_\varepsilon^4(1-|\varphi|^2)^2+ |{\rm
curl}\,A-H|^2\right)\md x.\nonumber
\end{equation}
\end{lem}

\section{Upper bound of the energy}\label{Section:UB}

\subsection{Main result}
The objective of this section is to establish the following upper
bound.

\begin{prop}\label{prop-upperbound}
Assume that
$\displaystyle\lim_{\varepsilon\to0}\frac{H}{|\ln\varepsilon|}=\lambda$
with $\lambda>0$. Given $\mu\in\mathcal M(\Omega)\cap
H^{-1}(\Omega)$, there exists a family of configurations
$\{(\varphi_\varepsilon,A_\varepsilon)\}_\varepsilon$ in
$H^1(\Omega;\mathbb C)\times H^1(\Omega;\mathbb R^2)$ such that
$\|\varphi_\varepsilon\|_{L^\infty(\Omega)}\leq 1$,
\begin{equation}\label{cv-mesure}
\frac{\mu(\varphi_\varepsilon,A_\varepsilon)}{H}\rightharpoonup\mu\quad{\rm
in}~\left(C_0^{0,\gamma}(\Omega)\right)^*\quad\forall~\gamma\in(0,1),
\end{equation}
and
$$\limsup_{\varepsilon\to 0} F_{\varepsilon,H}(\varphi_\varepsilon,A_\varepsilon)\leq\frac1\lambda\int_{\Omega}p(x)|\mu|\,\md
x+\int_\Omega\left(\frac1{p(x)}|\nabla
h_\mu|^2+|h_\mu-1|^2\right)\,\md x\,.$$
\end{prop}

\subsection{Analysis of a Green's potential}
This section is devoted to an analysis of  a Green's kernel, i.e.
a fundamental solution of the differential operator\break $-{\rm
  div}\,\left(\displaystyle\frac1{u_\varepsilon^2(x)}\nabla\right)+1$.
The existence and the properties of this function, taken from
\cite{AfSaSe, Stamp}, are given in the next lemma.

\begin{lem}\label{Green}
For every $y\in\Omega$ and $\varepsilon\in]0,1]$, there exists a
unique symmetric function $\Omega\times\Omega\ni(x,y)\mapsto
G_\varepsilon(x,y)\in\mathbb R_+$ such that~:
\begin{equation}\label{equation-Green}
\left\{\begin{array}{rl} -{\rm
div}\,\left(\displaystyle\frac1{u_\varepsilon^2(x)}\nabla_x
G_\varepsilon(x,y)\right)+G_\varepsilon(x,y)=\delta_y(x)&{\rm
in}~\mathcal
D'(\Omega),\\
G_\varepsilon(x,y)\big{|}_{x\in\partial\Omega}=0.
\end{array}\right.
\end{equation}
Moreover, $G_\varepsilon$ satisfies the following properties:
\begin{enumerate}
\item There exists a constant $C>0$ such that
$$G_\varepsilon(x,y)\leq C\left(|\,\ln|x-y|\,|+1\right)\,,\quad
\forall~(x,y)\in\overline\Omega\times\overline\Omega\setminus\Delta\,,~
\forall~\varepsilon\in]0,1]\,,$$ where $\Delta$ denotes the
diagonal in $\mathbb R^2$. \item The function $v_\varepsilon:
\Omega\times\Omega\ni (x,y)\mapsto
G_\varepsilon(x,y)+\displaystyle\frac{u_\varepsilon^2(x)}{2\pi}\ln|x-y|$
is in the class $C^1(\Omega\times\Omega\,;\mathbb R)$\,. \item
Given $q\in[1,2[$, there exists a constant $C>0$ such that
$$\|v_\varepsilon(\cdot,y)\|_{W^{1,q}(\Omega)} \leq C,
\quad
\forall~y\in\overline\Omega,\quad\forall~\varepsilon\in]0,1]\,.$$
\end{enumerate}
\end{lem}

One feature of the function $G_\varepsilon$ is that we know its
homogenized limit $G_0$.

\begin{lem}
\label{conv-G0}
 Let $1<q<2$. The function
$G_\varepsilon$ converges
 weakly in
$W^{1,q}(\Omega\times\Omega)$ to $G_0$ the solution of
\begin{equation}\label{equation-Green0}
\left\{\begin{array}{rl} -{\rm
div}\,\left(\displaystyle\frac1{p(x)}\nabla_x
G_0(x,y)\right)+G_0(x,y)=\delta_y(x)&{\rm in}~\mathcal
D'(\Omega),\\
G_0(x,y)\big{|}_{x\in\partial\Omega}=0.
\end{array}\right.
\end{equation}
Moreover, $G_\e$ converges locally uniformly to $G_0$ in $\Omega\times
\Omega\backslash \Delta$ where $\Delta$ is the diagonal of $\mathbb{R}^2$.
\end{lem}
\begin{proof}
Since $G_\varepsilon$ and $G_0$ are symmetric, it is sufficient to
establish weak convergence of $G_\varepsilon(\cdot,y)$ for
$y\in\Omega$
arbitrarly fixed.\\
Lemma~\ref{Green} assures that $G_\varepsilon(\cdot,y)$ is bounded
in $W^{1,q}(\Omega)$ uniformly with respect to $y$ and
$\varepsilon$. Hence, $G_\varepsilon(\cdot,y)$ converges weakly to
a function $G_0(\cdot,y)$ in $W^{1,q}(\Omega)$. Let us prove that
$G_0$ solves (\ref{equation-Green0}).\\
Let $\varphi\in C_0^\infty(\Omega)$. It is sufficient to prove
that as $\varepsilon\to0$,
\begin{eqnarray*}
&&\int_\Omega\frac{1}{u_\varepsilon^2(x)}\nabla_x
G_\varepsilon(x,y) \cdot\nabla _x\varphi(x)\,\md x\to
\int_\Omega\frac{1}{p(x)}\nabla_x
G_0(x,y)\cdot\nabla _x\varphi(x)\,\md x\\
&&\int_\Omega G_\varepsilon(x,y) \varphi(x)\,\md x\to \int_\Omega
G_0(x,y)\,\varphi(x)\,\md x\,.
\end{eqnarray*}
The second convergence is actually immediate since
$W^{1,q}(\Omega)$ comapctly embedds in $L^q(\Omega)$ and
$G_\varepsilon(\cdot,y)$ converges weakly to
$G_0(\cdot,y)$ in $W^{1,q}(\Omega)$. In order to prove the first
convergence we notice that,
\begin{eqnarray*}
&&\int_\Omega\frac{1}{u_\varepsilon^2(x)}\nabla_x
G_\varepsilon(x,y) \cdot\nabla _x\varphi(x)\,\md x-
\int_\Omega\frac{1}{p(x)}\nabla_x
G_0(x,y)\cdot\nabla _x\varphi(x)\,\md x\\
&&\hskip0.5cm=\int_\Omega\left(\frac{1}{u_\varepsilon^2(x)}-\frac1{p(x)}\right)
\nabla_x G_\varepsilon(x,y) \cdot\nabla _x\varphi(x)\,\md x\\
&&\hskip1cm+
\int_\Omega\frac{1}{p(x)}\nabla_x \left(\mathcal
G_\varepsilon(x,y)-G_0(x,y)\right)\cdot\nabla _x\varphi(x)\,\md
x\,,
\end{eqnarray*}
so it suffices to prove that
$$\int_\Omega\left(\frac{1}{u_\varepsilon^2(x)}-\frac1{p(x)}\right)
\nabla_x G_\varepsilon(x,y) \cdot\nabla _x\varphi(x)\,\md x\to
0\quad{\rm as}~\varepsilon\to0\,.$$ Using the uniform bound on
$G_\varepsilon(\cdot,y)$ in $W^{1,q}$, we write,
\begin{eqnarray*}
&&\int_\Omega\left|\left(\frac{1}{u_\varepsilon^2(x)}-\frac1{p(x)}\right)
\nabla_x G_\varepsilon(x,y)
\cdot\nabla _x\varphi(x)\right|\,\md x\\
&&\leq \left(\int_\Omega
\left(\frac{1}{u_\varepsilon^2(x)}-\frac1{p(x)}\right)^p|\nabla\varphi|^p\,
\md x \right)^{1/p}\left(\int_\Omega|\nabla_x
G_\varepsilon(x,y)|^q\,\md x
\right)^{1/q}\\
&&\leq C\left(\int_\Omega
\left(\frac{1}{u_\varepsilon^2(x)}-\frac1{p(x)}\right)^p|\nabla\varphi|^p\,
\md x \right)^{1/p}\,,\quad\frac1p+\frac1q=1\,.
\end{eqnarray*}
Now, we use the asymptotic behavior of $u_\varepsilon^2$, which
converges pointwise to $p(x)$ and remains uniformly bounded from
above and below. Therefore, Lebesgue convergence theorem
yields the desired convergence.\\
To complete the proof, we refer to \cite[Lemma~3.1]{AfSaSe} to find  the uniform convergence of
$G_\e$ to $G_0$.
\end{proof}

Next we cite a result from \cite[Lemma~5.5]{Kach-v}.

\begin{lem}\label{integration}
Let $(f_\varepsilon)_{\varepsilon\in]0,1]} \subset C([0,1],\mathbb
R_+)$ be a family of continuous functions. Assume that  there
exists a constant $C>0$ such that
$$\|f_\varepsilon\|_{L^1([0,1])}\leq C\ln|\ln\varepsilon|\,,\quad\forall~\varepsilon\in
]0,1]\,.$$ There exist constants $K>0$ and $c_0\in]0,1[$ such
that, given a family $(N(\varepsilon))\subset\mathbb N$ satisfying
$N(\varepsilon)\gg1$,  there exists a family
$(\delta(\varepsilon))\subset ]0,1[$ and a sequence
$(t_{m}^\varepsilon)_{m\in\mathbb N}\subset ]0,1[$
and
\begin{eqnarray*}
&&|f_\varepsilon(t_i^\varepsilon)|\leq K\ln|\ln\varepsilon|\,,
\quad \frac{c_0}{N(\varepsilon)}
\leq\left|t_{i+1}^\varepsilon-t_i^\varepsilon\right|\leq
\delta(\varepsilon)+\frac{c_0}{N(\varepsilon)}\,,\\
&&\forall~i\in\{1,2,\cdots,N(\varepsilon)\}\,,\quad\forall~\varepsilon\in
]0,1]\,.
\end{eqnarray*}
\end{lem}

The next proposition provides us with  points enjoying useful
properties. These points will serve to be the centers of the
vortices of the test configuration that we shall construct in the
next section.

\begin{prop}\label{prop-vp}
Let $\mu\in C_0^0(\Omega)$ be  continuous and compactly supported
in $\Omega$. Assume that
$(n(\varepsilon))_{\varepsilon\in(0,1]}\subset \mathbb N$ is a
family of integers such that $c_1|\ln\varepsilon|\leq
n(\varepsilon)\leq c_2\varepsilon^{-2}$ for
constants $c_1$ and $c_2$ independent of $\varepsilon\in(0,1]$.\\
If the restriction $\tilde\mu$ of $\mu$ to $S_1$ (respectively
$S_2$) is different from $0$, then there exist  constants  $c>0$, $C>0$,
a family of integers $N(\varepsilon)\in\mathbb N$, points
$a_i^\varepsilon$ in $S_1$ (respectively $S_2$) and degrees
$d_i^\varepsilon=\pm1$ such that:
\begin{enumerate}
\item $N(\varepsilon)=n(\varepsilon)(1+o(1))$ as
$\varepsilon\to0$\,;
\item $|a_i^\varepsilon-a_j^\varepsilon|\geq
\displaystyle\frac{C}{\sqrt{n(\varepsilon)}}\geq c\varepsilon$ for all
  $i\not=j$\,;
\item
${\rm dist}\left(a_i^\varepsilon,\partial
S_1\right)\geq\displaystyle\frac{\ln|\ln\varepsilon|}{|\ln\varepsilon|}$
and
  ${\rm dist}(a_i^\varepsilon,\partial\Omega)>c$ for
  all $i$\,;
\item $|v_\varepsilon(a_i^\varepsilon,a_i^\varepsilon)|\ll
|\ln\varepsilon|$ for all $i$\,; \item  $\displaystyle
\frac1{N(\varepsilon)}\sum_{i=1}^{N(\varepsilon)}
\mu_i^\varepsilon\rightharpoonup
2\pi\frac{\tilde\mu}{\|\tilde\mu\|}$ in the weak sense of
measures, for $\mu_i^\varepsilon$ any measure supported in
$\overline{B(a_i^\varepsilon,c\varepsilon)}$, of constant sign,
and
  such that $\mu_i^\varepsilon(\Omega)=2\pi d_i^\varepsilon$.
\end{enumerate}
\end{prop}
\paragraph{\bf Proof.}
For $i=1,2$, we denote by $S_i^\varepsilon=\{x\in S_i~:~{\rm
dist}(x,\partial S_1)\geq
\frac{\ln|\ln\varepsilon|}{|\ln\varepsilon|}\}$. We only present the
construction in $S_1^\e$ as it holds exactly in $S_2^\e$.\\
We partition $S_1^\varepsilon$ into squares $K$ of side-length
$\ell(\varepsilon)$ where $\ell(\varepsilon)$ is chosen such that,
$$\frac1{\sqrt{n(\varepsilon)}}\ll\ell(\varepsilon)\ll1\quad{\rm as}~
\varepsilon\to0\,.$$ We denote by $\mathcal K(\varepsilon)$ the
family of such squares
that lie entirely inside $S_1^\varepsilon$.\\
For $K\in\mathcal K(\varepsilon)$, we set
$$\lambda_K=n(\varepsilon)\frac{|\mu(K)|}{M_\varepsilon}\,,$$
where
$$M_\varepsilon=\sum_{K\in\mathcal K(\varepsilon)}
|\mu(K)|\,.$$ Since $\mu$ is continuous with compact support in
$\Omega$, and $\ell(\varepsilon)$ tends to $0$ as
$\varepsilon\to0$, we get
that $M_\varepsilon\to \|\tilde\mu\|$ as $\varepsilon\to0$.\\
Now $\displaystyle\sum_{K\in\mathcal K(\varepsilon)}\lambda_K=n$,
hence
  Lemma~7.4
of \cite{SaSe} provides us with nonnegative integers
$m_\varepsilon(K)$ such that
\begin{equation}\label{eq-mK}
  \sum_{K\in\mathcal
  K(\varepsilon)}m_\varepsilon(K)=n(\varepsilon)\,,\end{equation}
and
\begin{equation}\label{eq-mK-est}
\left|m_\varepsilon(K)-n\frac{|\mu(K)|}{M_\varepsilon}\right|<1\,.
\end{equation}
Since $\mu$ is bounded, we may find $M>0$ such that $|\mu|(K)\leq
M\ell(\varepsilon)^2$ for all $K\in\mathcal K(\varepsilon)$.
Consequently, we infer from (\ref{eq-mK-est}) that
$$m_\varepsilon(K)\leq 1+nM\ell^2=\mathcal O(nM\ell^2)\,,$$
(recall that $\frac1{n}\ll\ell\ll1$ as $\varepsilon\to0$).\\
Now we decompose $\mathcal K(\varepsilon)=\mathcal
K_1(\varepsilon)\cup \mathcal K_2(\varepsilon)$, where
$K\in\mathcal K_2(\varepsilon)$ if and only if\break
${\rm dist}(\partial K,\partial S_1)\geq \frac R2$ where $R>0$ is picked small enough so that
each connected component of $S_1$ contains a disc of radius $R$.\\
Now assume that $K\in \mathcal K_1(\varepsilon)$. Let
$p_\varepsilon(K)$ be the least integer less than
$\sqrt{m_\varepsilon(K)}$. Then we may pick $p_\varepsilon(K)$
points $b_i^\varepsilon$ lying  on one side $\mathcal L$ of $K$
and such that
\begin{equation}\label{eq-bi}
|b_i^\varepsilon-b_j^\varepsilon|\geq
\frac{C\ell(\varepsilon)}{\sqrt{p_\varepsilon(K)}}\geq
  \frac{C}{\sqrt{n(\varepsilon)}}\,.\end{equation}
Next, we can construct $p_\varepsilon(K)$ segments
$L_i^\varepsilon$ in $K$, each passing by $b_i^\varepsilon$ and
perpendicular to the
side $\mathcal L$ (hence each $L_i$ has length $\ell(\varepsilon)$).\\
Since $K\in \mathcal K_1(\varepsilon)$, we may find closed smooth
curves $\gamma_i\supset L_i$ in $\Omega$ such that each $\gamma_i$
is the boundary of a simply connected set $U_i\subset\Omega$ and
$2\pi R\leq |\gamma_i|\leq |\partial S_1|$. By this way, referring
to \cite[Proof of (5.4)]{Kach-v}, we may get a constant $C'>0$
independent from $\varepsilon$ and such that, for all $i$,
$$\int_{\gamma_i}|v_\varepsilon(x,x)|\,\md x \leq
C'\ln|\ln\varepsilon|\,.$$ In particular, it holds that
$$\int_{L_i}|v_\varepsilon(x,x)|\,\md x\leq C'\ln|\ln\varepsilon|\,.$$
Let us parameterize $L_i$ by $(0,\ell(\varepsilon)]\ni s \mapsto
x(s)\in L_i$, and let us define the rescaled continuous function
$$(0,1]\ni t\mapsto f_\varepsilon(t)=
v_\varepsilon\left(x(\ell(\varepsilon)t),x(\ell(\varepsilon)t)\right)\,,$$
so that $\|f_\varepsilon\|_{L^1(0,1)}\leq
C''(\ell(\varepsilon))^{-1}\ln|\ln\varepsilon|$. Invoking the
result of Lemma~\ref{integration}, we may pick $p_\varepsilon(K)$
points $(b_{i,j}^\varepsilon)\subset L_i$ such that\,\footnote{We
apply the lemma with
  $N(\varepsilon)=\sqrt{n(\varepsilon)}$ so that we get
  $\sqrt{n(\varepsilon)}\geq p_\varepsilon(K)$ points.}
$$|v_\varepsilon(b_{i,j}^\varepsilon,b_{i,j}^\varepsilon)|\leq
K'\frac{\ln|\ln\varepsilon|}{\ell(\varepsilon)}$$ and
$$|b_{i,j}^\varepsilon-b_{i,k}^\varepsilon|\geq \frac{c_0}
{\ell(\varepsilon)\sqrt{n(\varepsilon)}}\geq
\frac{c_0}{\sqrt{n(\varepsilon)}} \,.$$ Denoting by
$(a_k^\varepsilon)$ the family of the  constructed points in all
the segments $L_i$, we have actually constructed
$[p_\varepsilon(K)]^2$ points $(a_k^\varepsilon)\subset K$ such
that
$$|v_\varepsilon(a_k^\varepsilon,a_k^\varepsilon)|\leq
K'\frac{\ln|\ln\varepsilon|}{\ell(\varepsilon)}\,,\quad
|a_k^\varepsilon-a_j^\varepsilon|\geq
\frac{C}{\sqrt{n(\varepsilon)}}\,.$$ Recalling the assumptions on
$n(\varepsilon)$ and $\ell(\varepsilon)$, we get actually the two
desired properties (2) and (3) stated in Proposition~\ref{prop-vp}
above, with the constant $c>0$ chosen
sufficiently small that $c<{\rm dist}({\rm supp}\,\mu\,,\partial\Omega)$.\\
It remains now to continue the construction of the points
$(a_i^\varepsilon)$ filling the squares  $K\in \mathcal
K_2(\varepsilon)$. Notice that Corollary~5.3 of \cite{Kach-v} (more
precisely an adjustment of its proof) provides us with a constant
$C>0$ such that, for all $K\in\mathcal K_2(\varepsilon)$ and
$\varepsilon\in(0,1)$, we have,
$$\|v_\varepsilon(x,x)\| _{L^\infty(K)}\leq C\,.$$
Hence, it is sufficient to construct any well separated
$m_\varepsilon(K)$
 points in this case. This is exactly the case of
 \cite[p.~144]{SaSe}. Now, the integer $N(\varepsilon)$ is defined as
$$N(\varepsilon)=\sum_{K\in\mathcal K_1(\varepsilon)}
\left(p_\varepsilon(K)\right)^2+\sum_{K\in\mathcal
K_2(\varepsilon)} m_\varepsilon(K)\,.$$ That
$N(\varepsilon)=n(\varepsilon)(1+o(1))$ is due to the assumption
we made on $\ell(\varepsilon)$,
$\sqrt{n(\varepsilon)}^{-1}\ll\ell(\varepsilon)
\ll 1$ as $\varepsilon\to0$.\\
Now, having constructed the family of points $(a_i^\varepsilon)$, we
see that  property  (3)   stated in Proposition~\ref{prop-vp} is
just due to our  construction of the points being in
$S_\varepsilon^1$, and our choice of the constant $c$ being so
small  that ${\rm dist}({\rm supp}\,\mu\,,\partial\Omega)>c$.\\
Now we define the family $(d_i)$. If $\mu(K)\geq 0$, we assign the
degree $d_i=1$ to each $a_i\in K$, otherwise we associate the
degree $d_i=-1$. The proof of the last property (5) in
Proposition~\ref{prop-vp} is exactly as
that given in \cite[p.~145]{SaSe}.\hfill$\Box$\\

Proposition~\ref{prop-vp} will be used in the following context.
Let $\mu\in C_0^0(\Omega)$ be continuous and compactly supported
in $\Omega$. Take positive integers $n_1(\varepsilon)$ and
$n_2(\varepsilon)$ such that
$$c_1|\ln\varepsilon|\leq n_1(\varepsilon)+n_2(\varepsilon)\leq
 c_2\varepsilon^{-2}\,,$$
for positive constants $c_1$ and $c_2$ independent from
 $\varepsilon$.\\
If the restrictions $\check\mu$ and $\hat\mu$ of $\mu$ to $S_1$
and $S_2$ respectively  are both different from $0$, then we get
$N_1(\varepsilon)\sim n_1(\varepsilon)$ points $(a_i^\varepsilon)$
in $S_1$ and
 $N_2(\varepsilon)\sim n_2(\varepsilon)$
points $(b_i^\varepsilon)$ in $S_2$ satisfying properties (1)-(5)
of  Proposition~\ref{prop-vp}. In particular, we set
$$x_i^\varepsilon=\left\{
\begin{array}{ll}
a_i^\varepsilon&\forall~i\in\{1,2,\cdots,
n_1(\varepsilon)\}\,,\\
b_i^\varepsilon&\forall ~i
\in\{n_1(\varepsilon)+1,n_2(\varepsilon)+2,\cdots,
n_1(\varepsilon)+n_2(\varepsilon)\}\,,
\end{array}\right.$$
and we define the measures $\mu_i^\varepsilon$ by
\begin{equation}\label{eq:measure1}
\mu_i^\varepsilon(x)= \left\{
\begin{array}{cl}
\displaystyle\frac{2d_i}{c^2\varepsilon^2}&{\rm if}~ x\in
B(x_i^\varepsilon,c\varepsilon)
 \\
0&{\rm otherwise}\,.
\end{array}\right.\end{equation}
If $\check\mu$ (respectively $\hat\mu$) is zero, we may still choose
the points $a_i^\varepsilon$ (respectively $b_i^\varepsilon$)
arbitrarily so that properties (1)-(4) of Proposition~\ref{prop-vp}
are valid, and we define the corresponding
measures $\mu_i^\varepsilon$ to be zero by pure convention.\\
With these notations, we get as an immediate consequence that
\begin{equation}\label{eq:measure}
\mu_\e:=\frac1{n(\varepsilon)}
\sum_{i=1}^{N_1(\varepsilon)+N_2(\varepsilon)} \mu_i^\varepsilon
\rightharpoonup 2\pi\frac{\mu}{\|\mu\|}\quad{\rm
in~}(C_0^{0,\gamma}(\Omega))^*\,,
\end{equation}
 where
\begin{equation}\label{eq:n}
n(\varepsilon)=\left\{
\begin{array}{cl}
N_1(\varepsilon)&{\rm if}\quad{\rm supp}\,\mu\subset\overline S_1\,,\\
N_2(\varepsilon)&{\rm if}\quad{\rm supp}\,\mu\subset\overline
S_2\,,\\
N_1(\varepsilon)+N_2(\varepsilon)&{\rm otherwise.}\end{array}
\right.
\end{equation}

\begin{lem}\label{corol-cont-green}
Under the hypotheses and notations above, one has
\begin{equation}
\limsup_{\varepsilon\to0}\left(\frac1{n(\varepsilon)}\right)^2\sum_{i\not=j}
\int_{B_i\times B_j} G_\varepsilon(x,y)\md \mu_i^\varepsilon(x)\md
\mu_j^\varepsilon(y)\leq
\frac{4\pi^2}{\|\mu\|^2}\int_{\Omega\times\Omega} G_0(x,y)\md
\mu(x)\md \mu(y)\,. \label{Ge-mu}
\end{equation}
Here, for all $i$, $B_i$ denotes the ball
$B(x_i^\varepsilon,c\,\varepsilon)$.
\end{lem}
\begin{proof}
Given  $\alpha>0$, let
$\Delta_{\alpha}=\{(x,y),\,\,\,|x-y|<\alpha\}$. We have the
following decomposition,
\begin{equation}\label{double-somme1}
\frac1{n^2}\sum_{i\not=j} \int_{B_i\times B_j} G_\varepsilon\,\md
\mu_i^\varepsilon\md
\mu_j^\varepsilon=\int_{\Omega\times\Omega\backslash\Delta_{\alpha}}
G_\e \md \mu_\e \md\mu_\e +\frac1{n^2}\sum_{i\not=j}
\int_{\Delta_\alpha} G_\varepsilon\,\md \mu_i^\varepsilon\md
\mu_j^\varepsilon\,.
\end{equation}
Knowing from Lemma~\ref{conv-G0} that $G_\e$ converges uniformly  to
$G_0$  in
$\Omega\times\Omega\backslash\Delta_{\alpha}$, we write using
(\ref{eq:measure}),
\begin{equation}
\label{cv-Ge}
\lim_{\e\rightarrow 0}
\int_{\Omega\times\Omega\backslash\Delta_{\alpha}} G_\e \md \mu_\e
\md\mu_\e= \frac{4\pi^2}{\|\mu\|^2}\int_{\Omega\times\Omega\backslash\Delta_{\alpha}} G_0 \md
\mu \,\md \mu\,.
\end{equation}
Let us estimate now the last term on the right hand side of
(\ref{double-somme1}). Since the supports of $\mu_i^\e$ and $\mu_j^\e$
are disjoint for $i\not=j$, we may write by Lemma~\ref{Green},
$$
\left|\frac1{n^2} \sum_{i\not=j} \int_{\Delta_{\alpha}} G_\e \md
\mu_i^\e \md\mu_j^\e\right|\leq \frac{C}{n^2} \sum_{i\not=j}
\int_{\Delta_{\alpha}} |\,\ln|x-y|\,|\,|\mu_i^\e|(x)
\,|\mu_j^\e|(y)\,\md x\,\md y\,,$$
where $C>0$ is a constant independent from $\alpha$ and $\e$.\\
Proposition~\ref{prop-vp} provides us that the points $a_i^\e$  are
well separated, i.e.  $|a_i^\e-a_j^\e|\geq
\displaystyle\frac{C}{\sqrt{n}}$. This actually permits us to write
(see \cite[p. 147]{SaSe} for details),
\begin{equation}\label{A-S-S}
\left|\frac1{n^2}
\sum_{i\not=j}
\int_{\Delta_{\alpha}} G_\e \md \mu_i^\e
\md\mu_j^\e\right|\leq
C\int_{\Delta_\alpha}\left(|\,\ln|x-y|\,|+1\right)
\md x\,\md y\,.
\end{equation}
Substituting (\ref{cv-Ge}) and (\ref{A-S-S}) in
(\ref{double-somme1}), we deduce that,
\begin{eqnarray*}
\limsup_{\e\to0}
\frac{1}{n^2}\sum_{i\not=j}\int_{B_i\times B_j} G_\e \md\mu_i^\e(x)
\md\mu_j^\e(y)&\leq&
\frac{4\pi^2}{\|\mu\|^2}
\int_{\Omega\times\Omega\backslash\Delta_{\alpha}} G_0 \md
\mu \,\md \mu\\
&&+C\int_{\Delta_\alpha}\left(|\,\ln|x-y|\,|+1\right)
\md x\,\md y\,,
\end{eqnarray*}
with $\alpha$ being arbitrary in the interval $(0,1)$.
Making $\alpha\to0$ (recall
that $\ln|x-y|$ is in $L^1$), we obtain
the desired bound of the lemma.
\end{proof}

\begin{definition}\label{def1}
A family of points $(a_i)\subset\Omega$ satisfying Properties
(2)-(4) stated in Proposition~\ref{prop-vp} is said to be a
well-distributed family.
\end{definition}

Now, let $h_\varepsilon:\Omega\longrightarrow]0,1[$ be the solution
of the equation:
\begin{equation}\label{London*}
-{\rm div}\left(\frac1{u_\varepsilon^2}\nabla
h_\varepsilon\right)+h_\varepsilon=0\quad{\rm in}~\Omega\,,\quad
h_\varepsilon=1\quad{\rm on}~\partial\Omega\,,
\end{equation}
where $u_\varepsilon$ is introduced in Theorem~\ref{V-thm-kach3}. We
define also,
\begin{equation}\label{eq-J0}
J_0(\varepsilon)=\int_{\Omega}\left(\frac{1}{u_\varepsilon^2}|\nabla
h_{\varepsilon}|^2+|h_\varepsilon-1|^2\right)\,\md x\,.
\end{equation}

Next, we state a remarkable  energy-splitting  due to
Bethuel-Rivi{\`e}re \cite{BeRi}. We find it in
\cite[Lemma~5.7]{Kach-v}.

\begin{lem}\label{Beth-Riv}
Consider $(\varphi,A)\in H^1(\Omega;\mathbb C)\times
H^1(\Omega;\mathbb R^2)$ and define
$$A'=A-\frac{H}{u_\varepsilon^2}\nabla^\bot h_\varepsilon.$$
Then we have the decomposition of the energy,
\begin{equation}
\mathcal F_{\varepsilon,H}(\varphi,A)=H^2
J_0(\varepsilon)+\mathcal L_{\varepsilon,H}(\varphi,A')+\mathcal
R_0 +2 H\int_\Omega (h_\varepsilon-1) \mu(\varphi,A'),
\label{eq-split}
\end{equation}
where
\begin{equation}\label{eq-R0}
\mathcal R_0=H^2\int_\Omega \frac{1}{u_\e^2}(|\varphi|^2-1) |\nabla
h_\varepsilon|^2,\quad\quad \mu(\varphi,A')=h'+ {\rm
  curl}(i\varphi\,,\,(\nabla-iA')\varphi).
\end{equation}
Here, the functional $\mathcal L_{\varepsilon,H}$ is
\begin{equation}\label{eq-L}
\mathcal L_{\varepsilon,H}(\varphi,A')=\int_\Omega
u_\varepsilon^2|(\nabla-iA')\varphi|^2+|{\rm
curl}\,A'|^2+\frac1{2\varepsilon^2}
u_\varepsilon^4(1-|\varphi|^2)^2.
\end{equation}
\end{lem}

\begin{lem}
Given $\varepsilon\in (0,1)$, a well-distributed
family\,\footnote{See Definition~\ref{def1} above.} of $n=n_1+n_2$
points $(a_i)\subset\Omega$ ($n_i$ points in $S^i_\e$) together with
degrees $(d_i)\subset\{-1,1\}$, there exists a configuration
$(\varphi, A)$ such that, $\mu_i^\varepsilon$ being the uniform
measure on $\partial B_i=\partial B(a_i, c\varepsilon)$ of mass
$2\pi d_i$, and letting
 \[
\mu^\e=\sum_{i=1}^n \mu_i^\e,
\]
we have for some $\alpha\in [0,1[$,
\begin{equation}\label{autreformule}
\begin{split}
\mathcal{L}_{\e, H}(\varphi,A')& \leq c
 \e^{\alpha} n |\ln\e|^2+c' n \,\,o(|\ln\e|)+2\pi |\ln\e|
\sum_{i=1}^{n}p(a_i)\\
 &+\sum_{i\not=j} \int\int G_\varepsilon(x,y)\md \mu_i^\e(x)\md
\mu_j^\e(y)+\mathcal{O}(n),
\end{split}
\end{equation}
 \begin{equation}\label{formule-F}
 \begin{split}
\mathcal{F}_{\e, H}(\varphi,A)&=\mathcal{L}_{\e,
H}(\varphi,A')+H^2 J_0(\e)+2H\int_\Omega (h_\e-1) \md \mu(x)\\
&+\mathcal{O}\left(n \e H+n\e^2 H^2+(n^{\frac{1}{2}} \e H+\e H^2)
\sqrt{\mathcal{L}_{\e, H}(\varphi,A')}\right).
\end{split}
\end{equation}
Moreover, for any $0<\gamma\leq 1$, it holds that,
\begin{equation} \label{controlemesure} \frac{1}{n}
\|\mu(\varphi, A)-\mu\|_{C^{0,\gamma}(\Omega)}\leq C \e^{\gamma}
\left( 1+ \varepsilon H+ \sqrt{\frac{\mathcal{F}_{\e,H}(\varphi,
A')}{n}}\right).
\end{equation}
\label{lem-pv}
\end{lem}
\begin{proof} We construct a test configuration $(\varphi, A)$. We
define a function $h$ in $\Omega$ by $h=h'+H h_\e$ where
$h_\varepsilon$ has been introduced in (\ref{London*}) and $h'$ is
the solution of
\begin{equation}\label{V-Def-f'}
\left\{
\begin{array}{rl}
-{\rm div}\,\left(\displaystyle\frac1{u_\varepsilon^2}\nabla
h'\right)
+h'=\mu^\e&{\rm in}~\Omega,\\
h'=0&{\rm on}~\partial \Omega.\end{array}\right.\end{equation} Note
that $h'(x)=\int_{\Omega} G_\e(x,y) \md \mu^\e(y)$. As a
consequence, we have
\begin{equation} \label{formule-h'}
\int_\Omega \left(\frac{|\nabla h'|^2}{u_\e^2}+|h'|^2\right)\,\md
x=\int_\Omega\int_\Omega G_\e(x,y) \md \mu(x) \md \mu(y).
\end{equation}
Now we define an induced magnetic potential
$A=A'+\frac{H}{u_\varepsilon^2}\nabla^\bot h_\varepsilon$ by
  taking simply
$${\rm curl}\,A'=h'.$$
This choice is always possible  as one can take $A'=\nabla^\bot g$
with $g\in H^2(\Omega)$ such that
$\Delta g=h'$. We turn now to define an order
parameter $\varphi$ which we take in the form
\begin{equation}\label{V-vsol-psi'}
\varphi=\rho\,e^{i\phi},\end{equation} where $\rho$ is defined by
\begin{equation}\label{V-vsol-rho'}
\rho(x)=\left\{
\begin{array}{cll}
0&{\rm if}&x\in \cup_iB(a_i,c\varepsilon),\\
1&{\rm if}&x\not\in\cup_iB(a_i,2c\varepsilon),\\
\displaystyle\frac{|x-a_i|}\varepsilon-1&{\rm if}&\exists\,i~{\rm
s.t.}~x\in B(a_i,2c\varepsilon) \setminus B(a_i,c\varepsilon).
\end{array}\right.
\end{equation}
The phase $\phi$ is defined (modulo $2\pi$) by the relation:
\begin{equation}\label{V-vsol-phi'}
\nabla\phi-A'=-\frac{1}{u_\varepsilon^2}\nabla^\bot h'\quad{\rm
in}~ \Omega\setminus \cup_iB(a_i,c\varepsilon),
\end{equation}
and we emphasize here that we do not need to define $\phi$ in
regions where $\rho$ vanishes.\\
Having defined $(\varphi,A)$ as above, we estimate
$\mathcal{L}_{\e, H}(\varphi,A')$. Recall that
\begin{equation} \label{free energy1}
\mathcal{L}_{\e, H}(\varphi,A')=\int_{\Omega}
\left(u_\e^2|\nabla\rho|^2+\rho^2 u_\e^2
|\nabla\phi-A'|^2+|h'|^2+\frac{u_\e^4}{2\e^2}(1-\rho^2)^2\right)\md
x.
\end{equation}
From (\ref{V-vsol-rho'}) and using the uniform upper bound of
$u_\e$, it follows easily that,
\begin{equation}\label{controle-rho}
\int_{\Omega}\left(u_\e^2
|\nabla\rho|^2+\frac{u_\e^4}{2\e^2}(1-\rho^2)^2\right)\md x\leq C n.
\end{equation}
Thanks to (\ref{V-vsol-rho'}), (\ref{V-vsol-phi'}) and the
definition of $h=h'+Hh_\varepsilon$,  we have,
$$
u_\e^2 \rho^2 |\nabla\phi-A'|^2\leq u_\e^2 |\nabla\varphi
-A'|^2=\frac{|\nabla(h-H h_\e)|^2}{u_\e^2}\quad {\rm in}~
\Omega\setminus \cup_iB(a_i,c\varepsilon).\,,
$$
Replacing this in (\ref{free energy1}) and invoking
(\ref{controle-rho}), we find,
 \begin{equation}\label{free energy2}
\mathcal{L}_{\e, H}(\varphi,A')\leq \int_\Omega
\left(\frac{|\nabla(h-H h_\e)|^2}{u_\e^2}+|h-H h_\e|^2\right)\,\md
x+\mathcal{O}(n).
\end{equation}
Thanks to (\ref{formule-h'}), we deduce the  upper bound,
\begin{equation}\label{free energy3}
\mathcal{L}_{\e, H}(\varphi,A')\leq \int_{\Omega\times\Omega} G_\e
d\mu(x) d\mu(y)+\mathcal{O}(n)\,.
\end{equation}
We now decompose the double integral in two,
\begin{equation}\label{double-somme}
\int_{\Omega\times\Omega} G_\e d\mu^\e(x) d\mu^\e(y)=\sum_{i=1}^n
\int_{B_i\times B_i} G_\e \md\mu_i^\e(x)
\md\mu_i^\e(y)+\sum_{i\not=j}\int_{B_i\times B_j} G_\e
\md\mu_i^\e(x) \md\mu_j^\e(y).
\end{equation}
Let us estimate the first term in the right hand side. Writing
$G_\e=v_\e-\frac{u_\e^2}{2\pi} \ln |x-y|$, we have
\begin{equation}
\int_{B_i\times B_i} G_\e(x,y) \md\mu_i^\e(x)
\md\mu_i^\e(y)=\int_{B_i\times B_i}
\left(v_\e(x,y)-\frac{u_\e^2(x)}{2\pi} \ln |x-y|\right)
\md\mu_i^\e(x) \md\mu_i^\e(y),
\end{equation}
Assuming $(x,y)\in {\rm supp}\, \mu_i\times {\rm supp}\, \mu_i$, we
get
\begin{eqnarray*}
&&\hskip-1cm\int_{B_i\times B_i} u_\e^2(x) \ln |x-y| \md\mu_i^\e(x)
\md\mu_i^\e(y)\\
&&\hskip1cm=\int_0^{2\pi}\int_0^{2\pi} u_\e^2(a_i+c\e e^{i\theta_1})
\ln|c\e e^{i\theta_1}-c\e e^{i\theta_2}| \md\theta_1
\md\theta_2\\
&&\hskip1cm=c_\e+2\pi \ln\e\int_0^{2\pi}u_\e^2(a_i+c\e
e^{i\theta_1}) d\theta_1\,.
\end{eqnarray*}
Here and in the sequel,  $c_\e$ or $C_\e$ denote constants bounded
(uniformly in $\e$) from below and  above, and that may change from
one line to another.\\
The points $a_i$ being away from the boundary of $S_1$,
Lemma~\ref{interior} gives that\break $u_\varepsilon^2(a_i+c\e
e^{i\theta_1})$ is exponentially close to $p(a_i)$.  Thus, we may
write, upon considering the summation,
\begin{equation}
\label{sum}
\sum_{i=1}^n
 \int_{B_i\times B_i} \frac{u_\e^2(x)}{2\pi} \ln |x-y|
\md\mu_i^\e(x) \md\mu_i^\e(y)=C_\e n+\ln\e \sum_{i=1}^{n}p(a_i)\,.
\end{equation}
Let us now estimate $\sum_{i=1}^n \displaystyle \int_{B_i\times B_i}
|v_\e(x,y)| \md\mu_i^\e(x) \md\mu_i^\e(y)$. Referring to
Corollary~5.3 in
 \cite{Kach-v}, we know that for some $\alpha\in(0,1)$ and all $\eta\in(0,1)$,
 $$\|v_\e(\cdot,y)\|_{C^{0,\alpha}(\{x\in\Omega~:~{\rm dist}(x,\partial S_1)\geq \eta
 \}}\leq \frac{C_\alpha}{\eta^2}\,.$$
Consequently, the following estimate holds,
\begin{equation}
\begin{split}
\sum_{i=1}^n
 \int_{B_i\times B_i} |v_\e(x,y)| \md\mu_i^\e(x)
\md\mu_i^\e(y)&\leq \sum_{i=1}^n
 \int_{B_i\times B_i} \left(|v_\e(a_i,a_i)|
+C \frac{\e^{\alpha}}{\eta^2}\right)\md\mu_i^\e(x)
 \md\mu_i^\e(y)\\
 &= 4\pi^2  C n \frac{\e^{\alpha}}{\eta^2}+4\pi^2 n
 |v_\e(a_i,a_i)|,
\end{split}
\end{equation}
with
$\eta=\frac{2\ln|\ln\varepsilon|}{|\ln\varepsilon|}$.\\
By our hypotheses, we know that $|v_\e(a_i,a_i)|\ll|\ln\e|$. Thus,
\begin{equation}\label{cont-v}
\sum_{i=1}^n
 \int_{B_i\times B_i} |v_\e(x,y)| \md\mu_i^\e(x) \md\mu_i^\e(y)\leq c
 \e^{\alpha} n |\ln\e|^2+c' n \,\,o(|\ln\e|).
 \end{equation}
Combining (\ref{sum}) together with (\ref{cont-v}), we get
(\ref{autreformule}).\\
 The
proof of the  properties (\ref{formule-F})- (\ref{controlemesure})
are exactly as that given in \cite[p.~140-142]{SaSe}.\end{proof}

\subsection{Proof of proposition \ref{prop-upperbound}, completed}
We first assume $\mu\neq 0$ is a continuous and compactly supported
function. Let $n_1=[\frac{H}{2\pi} |\mu|(S_1)]$ and
$n_2=[\frac{H}{2\pi} |\mu|(S_2)]$ where $[\cdot]$ denotes the
integer part. We take $n=n_1+n_2$. Since $S_1$ and $S_2$ are
disjoint and cover $\Omega$, we have
$\|\mu\|=|\mu|(\Omega)=|\mu|(S_1)+|\mu|(S_2)$ and consequently,
\begin{equation}
\label{approx-n}
 \frac{n}{H}\backsim \frac{\|\mu\|}{2\pi}.
\end{equation}
Proposition~\ref{prop-vp} provides us with a well-distributed family
of points that serves as an input in Lemma~\ref{lem-pv}. Thus, we
get configurations $(\varphi_\e,A_\e)$ and associated measures
$\mu_i^\e$ such that
\[
\mu_\e:=\frac{1}{H}\sum_{i=1}^{n_1+n_2} \mu_i^{\e}\rightharpoonup
\mu\,, \quad{\rm weakly~in~}\left(C_0^{0,\gamma}\right)^*\,, \] and
the estimates (\ref{Ge-mu}), (\ref{autreformule})-(\ref{formule-F})
and (\ref{controlemesure}) being all valid.\\
Since $H$ has the order of $|\ln\e|$, Lemma~\ref{corol-cont-green}
and estimates (\ref{autreformule})-(\ref{approx-n}) yield,
$$
\mathcal{L}_{\e, H}(\varphi,A')\leq 2\pi|\ln\e|
\sum_{i=1}^{n}p(a_i)+4\pi^2\frac{n^2}{\|\mu\|^2}\int_{\Omega}\int_{\Omega}
G_0(x,y)\md \mu(x)\md \mu(y)+o(H^2),
$$
Inserting the particular choices of $n_1$, $n_2$, $n=n_1+n_2$ and
using the definition of $p$ being piecewise constant,we get,
 \begin{equation}\label{autreformule3}
\mathcal{L}_{\e, H}(\varphi,A')\leq H |\ln\e|\int_{\Omega} p(x)
|\mu|(x) \md x+H^2\int_{\Omega}\int_{\Omega} G_0(x,y)\md \mu(x)\md
\mu(y)+o(H^2).
\end{equation}
Again the hypothesis on $H=\mathcal{O}(|\ln\e|)$ gives that the
remainder terms in (\ref{formule-F}) are $o(1)$, leading thus to
\begin{eqnarray*}
{\mathcal F}_{\e,H}(\varphi,A)&\leq& H^2 J_0(\e)+H
|\ln\e|\int_{\Omega} p(x)
|\mu|(x) \md x\\
&&
+2 H^2 \int_\Omega (h_\e-1) \md
\mu_\e+H^2\int_{\Omega}\int_{\Omega} G_0(x,y)\md \mu(x)\md
\mu(y)+o(H^2).\nonumber
\end{eqnarray*}
Using (\ref{eq-J0}), we rewrite the preceding formula in the
following explicit form,
\begin{equation}\label{autreformule4}
\begin{split}
{\mathcal F}_{\e,H}(\varphi,A)& \leq  H^2 \int_\Omega
\left(\frac{|\nabla h_\e|^2}{u_\e^2}+|h_\e-1|^2\right)\,\md x +H
|\ln\e|\int_{\Omega} p(x) |\mu|(x) \md x\\
& +2 H^2 \int_\Omega (h_\e-1) \md \mu+H^2\int_{\Omega\times\Omega}
G_\e(x,y)\md \mu(x)\md
\mu(y)\\
&+2 H^2 \int_\Omega (h_\e-1) (\md \mu_\e-\md \mu)+
H^2\int_{\Omega\times\Omega} (G_0-G_\e) \md \mu \md \mu+o(H^2).
\end{split}
\end{equation}
Let us define  the functions $U_{\mu,\e}$ and $h_{\mu,\e}$ by
\begin{equation}\label{Umu-e}
\left\{
\begin{array}{rl}
-{\rm div}\,\left(\displaystyle\frac1{u_\e^2}\nabla
U_{\mu,\e}\right)
+U_{\mu,\e}=\mu&{\rm in}~\Omega,\\
U_{\mu,\e}=0&{\rm on}~\partial
\Omega.\end{array}\right.
\end{equation}
\begin{equation}\label{hmu-e}
\left\{
\begin{array}{rl}
-{\rm div}\,\left(\displaystyle\frac1{u_\e^2}\nabla
h_{\mu,\e}\right)
+h_{\mu,\e}=\mu&{\rm in}~\Omega,\\
h_{\mu,\e}=1&{\rm on}~\partial \Omega.\end{array}\right.
\end{equation}
Remarking that $U_{\mu,\e}(x)=\displaystyle\int_{\Omega} G_\e(x,y)
\md \mu(y)$, we get as a consequence
\begin{equation}\label{formule-hmu-e}
\int_\Omega \left(\frac{|\nabla
U_{\mu,\e}|^2}{u_\e^2}+|U_{\mu,\e}|^2\right)\,\md
x=\int_\Omega\int_\Omega G_\e(x,y) \md \mu(x) \md \mu(y)\,.
\end{equation} Writing $h_{\mu,\e}-1=U_{\mu,\e}+h_\e-1$ and
replacing (\ref{formule-hmu-e}) in
in (\ref{autreformule4}) leads, after some calculations  to
\begin{equation}\label{autreformule6}
\begin{split}
\frac{{\mathcal F}_{\e,H}(\varphi,A)}{H^2}& \leq \int_\Omega
\left(\frac{|\nabla
h_{\mu,\e}|^2}{u_\e^2}+|h_{\mu,\e}-1|^2\right)\md x+
\frac{|\ln\e|}{H} \int_\Omega p(x)\, |\mu|(x) \md x\\
& +2 \int_\Omega (h_\e-1) (\md \mu_\e-\md
\mu)+\int_{\Omega\times\Omega} (G_0-G_\e) \md \mu \md \mu+o(1).
\end{split}
\end{equation}
We define $h_{\mu}$ by
\begin{equation}\label{hmu}
\left\{
\begin{array}{rl}
-{\rm div}\,\left(\displaystyle\frac1{p(x)}\nabla h_{\mu}\right)
+h_{\mu}=\mu&{\rm in}~\Omega,\\
h_{\mu}=1&{\rm on}~\partial \Omega.\end{array}\right.\end{equation}
By a standard compactness argument similar to that given to
Lemma~\ref{conv-G0}, we check that,
\begin{equation}\label{eq:conv1}
\lim_{\e\longrightarrow 0}\int_\Omega \left(\frac{|\nabla
h_{\mu,\e}|^2}{u_\e^2}+|h_{\mu,\e}-1|^2\right)\md x=\int_\Omega
\left(\frac{|\nabla h_{\mu}|^2}{p(x)}+|h_{\mu}-1|^2\right)\md x.
\end{equation}
Using a Cauchy-Schwarz inequality, it follows  from Lemma
\ref{conv-G0}  and the boundedness of $\mu$ that
\begin{equation}\label{eq:conv2}
\lim_{\e\longrightarrow 0}\int_{\Omega\times\Omega} (G_0-G_\e) \md
\mu \md \mu=0.
\end{equation}
Next, noticing that $\|h_\e-1\|_{C^0_0(\Omega)}\leq 1$, it follows
immediately from the convergence of $\mu_\e$ to $\mu$ in
$(C_0^0(\Omega))^*$,
\begin{equation}
\lim_{\e\longrightarrow 0}\int_\Omega (h_\e-1) (\md \mu_\e-\md
\mu)=0\,. \label{cv-h-mu}
\end{equation}
Replacing  (\ref{eq:conv1})-(\ref{cv-h-mu}) in (\ref{autreformule6})
then using $\displaystyle\lim_{\e\to0}\frac{H}{|\ln\e|}= \lambda$,
we get finally,
\[
\limsup_{\e\longrightarrow 0}\frac{{\mathcal
F}_{\e,H}(\varphi,A)}{H^2}\leq \frac{1}{\lambda}\int_\Omega
p(x)|\mu|(x) \,\md x+\int_\Omega\left( \frac{1}{p(x)}|\nabla
h_{\mu}|^2+|h_{\mu}-1|^2\right) \md x .
\]
Moreover, since (\ref{controlemesure}) holds and ${\mathcal
L}_{\e,H}(\varphi,A')\leq C |\ln\e|^2$, we have
\[
\frac{1}{H}
\|\mu(\varphi,A)-\mu_\e\|_{(C_0^{0,\gamma}(\Omega))^*}\leq C
\e^{\gamma}\left(1+\sqrt{\frac{{\mathcal
L}_{\e,H}(\varphi,A')}{n}}\right)\leq o(1).
\]
We conclude that (\ref{cv-mesure}) holds, which finishes the proof
in the case where $\mu$ is a continuous and compactly supported
function. The general case where  $\mu\in \mathcal{M}(\Omega)\cap
H^{-1}(\Omega)$ follows by a standard  approximation argument, see
\cite[p.~149]{SaSe} for details.

\section{Lower Bound}\label{Sec:LB}

\subsection{Main result}
The objective of this section is to prove the lower bound stated
in Proposition~\ref{prop-lowerbound} below.\\
Given a family of configurations
$\{(\varphi_\varepsilon,A_\varepsilon)\}$, we denote by
\begin{equation}\label{eq:current}
j_\varepsilon=\left(i\varphi_\varepsilon,
(\nabla-iA_\varepsilon)\varphi_\varepsilon\right)\,,\quad
h_\varepsilon={\rm curl}\,A_\varepsilon\,.
\end{equation}

\begin{prop}\label{prop-lowerbound}
Assume that
$\displaystyle\lim_{\varepsilon\to0}\frac{H}{|\ln\varepsilon|}=\lambda$
with $\lambda>0$. Let $\{(\varphi_{\varepsilon},A_{\varepsilon})\}_n$ be
a family of configurations satisfying
$\mathcal F_{\varepsilon,H}(\varphi_\varepsilon,A_\varepsilon)\leq
C H^2$ and $\|\varphi_\varepsilon\|_{L^\infty(\Omega)}\leq 1$ for a
given constant $C>0$.\\
Then, up to the extraction of a subsequence $\varepsilon_n$
converging to $0$, one has,
$$\forall~\gamma\in(0,1)\quad
\frac{\mu(\varphi_{\varepsilon_n},A_{\varepsilon_n})}{H}\to\mu \quad{\rm
  in~}
\left(C^{0,\gamma}(\Omega)\right)^*\,,$$
$$\frac{j_{\varepsilon_n}}{H}\rightharpoonup j\,,\quad
\frac{h_{\varepsilon_n}}{H}\rightharpoonup h\quad
{\rm weakly~in~}L^2(\Omega)\,.$$
Moreover, $\mu={\rm curl}\,j+h$  and
\begin{equation}\label{eq:lb}
\liminf_{\varepsilon\to0} \frac{\mathcal
F_{\varepsilon,H}(\varphi_\varepsilon,A_\e)}{H^2} \geq
E_\lambda(\mu)+\int_\Omega\left( p(x)\left|j+\frac1{p(x)}\nabla^\bot
h_\mu\right|^2 +|h-h_\mu|^2 \right)\,\md x\,.
\end{equation}
Here, the energy  $E_\lambda$ and the function $h_\mu$ are introduced
in
(\ref{eq-E-lambda}) and (\ref{eq-h-mu}) respectively.
\end{prop}

\subsection{Vortex-balls}
In this section we construct suitable `vortex-balls' providing a
lower bound of the energy of minimizers of (\ref{V-EGL}). Recall
the decomposition of the energy in Lemma~\ref{V-lem-psi<u}, which
permits us to work with the `reduced energy functional' $\mathcal
F_{\varepsilon,H}$.\\
Notice that, by using
$(u_\varepsilon,0)$ as a test configuration for the function
(\ref{V-EGL}), we deduce an upper  bound of the form~:
\begin{equation}\label{V-upperbound}
\mathcal F_{\varepsilon,H}(\varphi,A)\leq CH^2\,,
\end{equation}
where $\varphi=\psi/u_\varepsilon$, $(\psi,A)$ a minimizer of
(\ref{V-EGL}), and $C>0$ a positive constant.\\
We recall the hypothesis that there exists a
positive constant $C>0$ such that the applied magnetic field $H$
satisfies
\begin{equation}\label{hypothesis-H}
H\leq C|\ln\varepsilon|\,.\end{equation}

The upper bound (\ref{V-upperbound}) provides us, as in
\cite{SaSe}, with the construction of suitable `vortex-balls'.

\begin{prop}\label{V-lem-vortexballs}
Assume the hypotheses (\ref{hypothesis-H}). Given an open set
$U\subset \Omega$ and a number $p\in]1,2[$, there exists a constant
$C>0$ and a finite family of disjoint balls $\{B((a_i,r_i)\}_{i\in
I}$ such that, $(\varphi,A)$ being a configuration  satisfying the
bound (\ref{V-upperbound}), the following properties hold:
\begin{enumerate}
\item $\overline{B(a_i,r_i)}\subset U$ for all $i$\,;
\item $w=\{x\in U~:~|\varphi(x)|\leq
  1-|\ln\varepsilon|^{-4}\}
\subset\displaystyle\displaystyle\bigcup_{i\in I}B(a_i,r_i)$.
\item $\displaystyle\sum_{i\in I}r_i\leq
C\,|\ln\varepsilon|^{-10}$. \item Letting $d_i$ be the degree of
the function $\varphi/|\varphi|$ restricted to $\partial
B(a_i,r_i)$ if $B(a_i,r_i)\subset\Omega$ and $d_i=0$
  otherwise, then  we have:
\begin{eqnarray}\label{V-Eq-LBest}
&&\hskip-0.5cm \int_{B(a_i,r_i)\setminus\omega}
u_\varepsilon^2|(\nabla-iA)\varphi|^2\,\md
x+\int_{B(a_i,r_i)}|{\rm
  curl}\,A-H|^2\,\md x\geq\\
&&\hskip3.5cm 2 \pi|d_i|
\left(\min_{B(a_i,r_i)}u_\varepsilon^2\right)
\left(|\ln\varepsilon|-C\ln|\ln\varepsilon|\right). \nonumber
\end{eqnarray}
\item $\left\|2\pi\displaystyle\sum_{i\in I} d_i\delta_{a_i}-{\rm
    curl}\big{(}A+(i\varphi,\nabla_A\varphi)\big{)}
\right\|_{W^{-1,p}_0(U)}\leq C|\ln\varepsilon|^{-4}.$
\end{enumerate}
\end{prop}

We follow the usual terminology and call the balls constructed in
Proposition~\ref{V-lem-vortexballs} `vortex-balls'. The proof of
Proposition~\ref{V-lem-vortexballs} is actually a simple
consequence of the analysis of \cite{SaSe}.

\subsection{Proof of Proposition~\ref{prop-lowerbound}}
Let us consider smooth and open sets $U_1\subset S_1$ and $U_2\subset
S_2$, and let us denote by $U$ their union.\\
Applying Proposition~\ref{V-lem-vortexballs} in $U_1$ and $U_2$
respectively, we get families of balls $B(a_i,r_i)$ and degrees $d_i$
such that,
$$
\frac{\mathcal
F_{\varepsilon,H}(\varphi_\varepsilon,A_\varepsilon,V_\varepsilon)}{H^2}
\geq 2\pi \frac{|\ln\varepsilon|}{H}
\sum_i\left(\min_{\overline{B(a_i,r_i)}}u_\varepsilon^2\right)
\frac{|d_i|}{H} +o(1)\,,
$$
where $V_\varepsilon=\displaystyle\cup_i B(a_i,r_i)$.\\
Recall that $u_\varepsilon^2$ converges uniformly to the function $p$
in $U$. Thus, we may rewrite the above lower bound in the following
form,
\begin{equation}\label{eq:lb-V}
\frac{\mathcal
F_{\varepsilon,H}(\varphi_\varepsilon,A_\varepsilon,V_\varepsilon)}{H^2}
\geq 2\pi \frac{|\ln\varepsilon|}{H}
\left(\frac{\sum_{a_i\in U_1} |d_i|}{H}+a\frac{\sum_{a_i\in U_2}|d_i|}H\right)
+o(1)\,.
\end{equation}
Using the bound $\|\varphi_\varepsilon\|_{L^\infty(\Omega)}\leq1$, we
write $|(\nabla-iA_\varepsilon)\varphi_\varepsilon|^2\geq
|\varphi_\varepsilon|^2|(\nabla-iA_\varepsilon)|^2\geq
|j_\varepsilon|^2$. Consequently, we have,
\begin{equation}\label{eq:lb-U}
\frac{\mathcal
F_{\varepsilon,H}(\varphi_\varepsilon,A_\varepsilon,U\setminus
V_\varepsilon)}{H^2} \geq \int_{U\setminus
V_\varepsilon}\left(u_\varepsilon^2\left|\frac{j_\varepsilon}{H}\right|^2+
\left|\frac{h_\varepsilon}{H}-1\right|^2\right)\,\md x\,.
\end{equation}
With $U=\Omega$, we infer from the bound
$\mathcal F_{\varepsilon,H}(\varphi_\varepsilon,A_\varepsilon)\leq
CH^2$ that up to the extraction of a subsequence,
 $\displaystyle\frac{j_\varepsilon}{H}$ and
$\displaystyle\frac{h_\varepsilon}{H}$ respectively
converge weakly to $j$ and $h$ in
$L^2(\Omega)$.\\
Thanks again to the bound (\ref{V-upperbound}), we deduce that
$\displaystyle\frac{\sum_{a_i\in U_1} |d_i|}{H}$ and
$\displaystyle\frac{\sum_{a_i\in U_2} |d_i|}{H}$ are bounded. Hence
the measures $\displaystyle\frac{\sum_{a_i\in U_1}
  d_i\delta_{a_i}}{H}$ and $\displaystyle\frac{\sum_{a_i\in U_2}
  d_i\delta_{a_i}}{H}$ are
 weakly compact in the sense of measures, and thus,
up to
extraction of subsequences, they respectively converge
to measures $\mu_1$ and $\mu_2$ in $(C_0^0(U_1))^*$ and $(C_0^0(U_2))^*$.
Thanks to the last property of
Proposition~\ref{V-lem-vortexballs}, we get upon setting
$\mu={\rm curl}\,j+h$,
$$\mu_1=\mu_{|_{U_1}}\,,\quad\mu_2=\mu_{|_{U_2}}\,.$$
Therefore, combining (\ref{eq:lb-V})-(\ref{eq:lb-U}) and using the
uniform convergence of $u_\varepsilon^2$ to $p$ in $U$, we deduce
that,
\begin{eqnarray*}
\frac{\mathcal
F_{\varepsilon,H}(\varphi_\varepsilon,A_\varepsilon,U)}{H^2} &\geq&
\frac{\mathcal
F_{\varepsilon,H}(\varphi_\varepsilon,A_\varepsilon,V_\varepsilon)}{H^2}
+\frac{\mathcal
F_{\varepsilon,H}(\varphi_\varepsilon,A_\varepsilon,U\setminus
V_\varepsilon)}{H^2}\\
&
\geq&\frac1{\lambda}\left(|\mu|(U_1)+a|\mu|(U_2)\right)+\int_U\left(
p(x)|j|^2+|h-1|^2\right)\,\md
x+o(1).
\end{eqnarray*}
Invoking $\mathcal
F_{\varepsilon,H}(\varphi_\varepsilon,A_\varepsilon,\Omega)\geq\mathcal
F_{\varepsilon,H}(\varphi_\varepsilon,A_\varepsilon,U)$, we deduce
that,
\begin{equation}\label{eq:lb-1}
\liminf_{\varepsilon\to0} \frac{\mathcal
F_{\varepsilon,H}(\varphi_\varepsilon,A_\varepsilon,\Omega)}{H^2}
\geq
\frac1{\lambda}\left(|\mu|(U_1)+a|\mu|(U_2)\right)+\int_U\left(p(x)|j|^2+|h-1|^2\right)\,\md
x\,.
\end{equation}
The left hand side of (\ref{eq:lb-1}) being independent from $U$, and
$U_1$, $U_2$  being arbitrary subsets of $S_1$ and $S_2$, we conclude
that
\begin{equation}\label{eq:lb-2}
\liminf_{\varepsilon\to0} \frac{\mathcal
F_{\varepsilon,H}(\varphi_\varepsilon,A_\varepsilon,\Omega)}{H^2}
\geq \frac1{\lambda}\int_\Omega p(x)\md |\mu|
+\int_\Omega\left(p(x)|j|^2+|h-1|^2\right)\,\md x\,.\end{equation}
Now to conclude, we write
$$j=-\frac1{p(x)}\nabla^\bot h_\mu+\left(j+\frac1{p(x)}\nabla^\bot
h_\mu\right)\,,\quad
h=h_\mu+(h-h_\mu)\,.$$
Upon substitution in the right hand side of (\ref{eq:lb-2}) and using
in particular the remarkable identities
$${\rm curl}\left(j+\frac1{p(x)}\nabla^\bot h_\mu\right)
+h-h_\mu=0\quad{\rm in~}\Omega\,,\quad h_\mu-1=0\quad{\rm
  on~}\partial\Omega\,,
$$
we get the desired conclusion (\ref{eq:lb}).
\hfill 

\begin{rem}\label{rem:conc}
Combining the upper and lower bounds of
Propositions~\ref{prop-upperbound} and
Proposition~\ref{prop-lowerbound}, then by uniqueness of the
minimizer $\mu_*$  of $E_\lambda$ (see Section~\ref{Sec:Min} below),
it is evident that $\mu_*=\mu$ and $h=h_{\mu_*}$. Here $\mu$ and $h$
are given in Proposition~\ref{prop-lowerbound} above.
\end{rem}

\section{Minimization of the limiting energy}\label{Sec:Min}
As we explained in the introduction, by convexity and lower
semi-continuity, the limiting energy (\ref{eq-E-lambda}) admits a
unique minimizer $\mu_*$ which is expressed by means of the unique
minimizer $h_*$ of (\ref{eq:dual-en}) as follows,
\begin{equation}\label{eq:mu*=h*}
\mu_*=-{\rm div}\left(\frac1{p(x)} \nabla h_*\right)+h_*\,.
\end{equation}

Proceeding as in \cite{SaSe3, SaSe}, we may get an equivalent
characterization of $h_*$. We write $H^1_1(\Omega)$ for the space of
Sobolev functions $u$ such that $u-1\in H^1_0(\Omega)$.
\begin{prop}\label{prop-min}
The minimizer $u_*$ of
$$
\min_{\substack{u\in H_0^1(\Omega)\\
 -{\rm div}(\frac{1}{p(x)}\nabla
  u)+u\in\mathcal{M}(\Omega)}}\int_\Omega\left( \frac{
p(x)}{\lambda}\left|-{\rm div}\left(\frac{1}{p(x)}\nabla
  u\right)+u+1\right|+\frac{|\nabla u|^2}{p(x)}+|u|^2\right)\,\md
  x\,,
$$
is also the unique minimizer of the dual problem
$$
\min_{ \substack{v\in H_0^1(\Omega)\\|v|\leq \frac{ p}{2\lambda}}}
\int_\Omega \left(\frac{|\nabla v|^2}{p(x)}+|v|^2+2 v\right)\,\md
x\,.
$$
For instance, $h_*=u_*+1$ minimizes the energy,
$$
\min_{\substack{f\in H_1^1(\Omega)\\(f-1)\geq -\frac{ p}{2\lambda}}}
\left(\int_\Omega \frac{|\nabla f|^2}{p(x)}+|f|^2\right)\,,
$$
and  satisfies $-{\rm div}\left(\frac{1}{p(x)}\nabla
  h_*\right)+h_*\geq 0$.
\end{prop}
\begin{proof}
Let us define the lower semi-continuous and convex functional
$$ \Phi(u)=
\int_\Omega \frac{p(x)}{2\lambda} \left|-{\rm div}
\left(\frac1{p(x)}\nabla u\right)+u+1\right|\,\md x$$ in the
Hilbert space $H=H^1_0(\Omega)$ endowed with the scalar product
$\langle f,g\rangle_H=\int_\Omega \frac{1}{p(x)}\nabla f\cdot \nabla
g+ f g$. Let us compute its conjugate $\Phi^*$, i.e.
$$\Phi^*(f)=\sup_{~\{g~:~\Phi(g)<\infty\}}\langle f,g\rangle-\Phi(g)\,.$$
Indeed, we have,
$$
\Phi^*(f) \geq  \sup_{\eta\in L^2}\int_\Omega f\eta\,\md
x-\frac{1}{2\lambda}\int_\Omega p(x) |\eta|\,\md x-\int_\Omega
f\,\md x\,,
$$from which we deduce that
\[
 \Phi^*(f)=\left\{
\begin{array}{ll}
-\displaystyle\int_\Omega f\,\md x&{\rm if}~|f|\leq \frac{ p}{2\lambda}\,,\\
+\infty&{\rm otherwise.}
\end{array}\right.
\]
By convex duality (see \cite[Lemma~7.2]{SaSe}),
$$\min_{u\in H}\left(\|u\|_H^2+2\Phi(u)\right)=-\min_{f\in
H}\left(\|f\|_H^2+2\Phi^*(-f)\right)\,,
$$ and minimizers  coincide. That the measure $\mu_*=-{\rm
div}\left(\frac{1}{p(x)}\nabla h_*\right) +h_*$ is positive is
actually a consequence of the weak maximum principle, see
\cite[p.~131]{Ka}. One may also follow step by step the proof given
in \cite{SaSe3}.
\end{proof}

Therefore, the limiting vorticity measure $\mu_*$ is positive, and
following \cite{SaSe3}, it can be expressed by means of the
coincidence set $w_\lambda=w_\lambda^1\cup w_\lambda^2$,
\begin{equation}\label{eq:w-l}
w_\lambda^1=\{ x\in \overline S_1~:~
1-h_*(x)=\frac{1}{2\lambda}\}\,,\quad w_\lambda^2=\{x\in\overline
S_2~:~\frac{1-h_*(x)}{a}=\frac1{2\lambda}\}\,,
\end{equation}
as follows,
\begin{equation}\label{eq:mu*}
\mu_*=\left(1-\frac{p(x)}{2\lambda}\right){\bf{1}}_{w_\lambda} \md
x\,,
\end{equation}
where $\mathbf 1_{w_\lambda}\md x$ denotes the Lebesgue measure
restricted to $w_\lambda$. Furthermore, $h_*$ (the minimizer of
(\ref{eq:dual-en})) solves,
\begin{equation}\label{eq:h*}
\left\{
\begin{array}{ll} -{\rm div}\left(\frac1{p(x)}\nabla
h_*\right)+h_*=0&{\rm in}~\Omega\setminus w_\lambda\\
h_*=1-\frac{p}{2\lambda}&{\rm in ~}w_\lambda\\
h_*=1&{\rm on~}\partial \Omega\,,
\end{array}\right.\end{equation}
and the regularity of $h_*$ is complicated, as $w_\lambda$ may
intersect the interior boundary $\partial S_1$ where $h_*$ has
gradient singularities (it is expected to satisfy a transmission
condition, see the radial case below). But, we know that in the
interior of  $S_1$ and $S_2$, $h_*$ is locally  $C^{1,\alpha}$ for
some exponent $\alpha\in(0,1)$. Remark for instance that the measure
$\mu_*$ is no more uniform and may be discontinuous in light of
$\mu_*=\left(1-\frac{1}{2\lambda}\right)\mathbf \md x$ in
$w_\lambda^1$ and $\mu_*=\left(1-\frac{a}{2\lambda}\right)\md x$ in
$w_\lambda^2$.\\
For the sake of a better understanding of the sets $w_\lambda$,
$w_\lambda^1$ and $w_\lambda^2$, we introduce the following critical
constants (we emphasize their dependence on the material parameter
$a$),
\begin{equation}\label{eq:l0}
\lambda_i(a)=\frac{1}{2\displaystyle \max_{x\in \overline{S_i}}
\left(\frac{1-h_0(x)}{p(x)}\right)},\quad\forall~i\in\{1,2\}\,,\quad
\lambda_0(a)=\min(\lambda_1(a),\lambda_2(a))\,.
\end{equation}
Here we recall  that  $h_0$ is the solution of $-{\rm
div}\left(\frac{1}{p(x)}\nabla
  h_0\right)+h_{0}=0$ in $\Omega$ and $h_0=1$ on $\partial\Omega$.
The maximum principle gives that $0<h_0<1$ in $\Omega$. We denote
also by $w_0=w_0^1\cup w_0^2$
the set $w_{\lambda_0(a)}=w_{\lambda_0(a)}^1\cup w_{\lambda_0(a)}^2$ introduced in (\ref{eq:w-l}).\\
The following proposition follows from a weak maximum principle
\cite[p.~131]{Ka}, see \cite{SaSe3} for a detailed proof (modulo
necessary adjustments).

\begin{prop}\label{prop:w-l}
 \begin{enumerate}
\item $w_\lambda$ is increasing with respect to $\lambda$ and
$\cup_{\lambda>0}w_\lambda=\Omega$\,;
\item $w_\lambda^1$ and $w_\lambda^2$ are disjoint\,;
\item If
$\lambda<\lambda_0(a)$ then $h_*=h_0$, $\mu_*=0$ and
$w_\lambda=\emptyset$\,;
\item If $\lambda=\lambda_0(a)$ then $h_*=h_0$, $\mu_*=0$ and
$w_\lambda=w_0$\,;
\item If $\lambda>\lambda_0(a)$ then $\mu_*\not=0$\,;
\item If $\lambda<\lambda_i(a)$ for some $i\in\{1,2\}$, then $w_{\lambda}^i=\emptyset$.
\end{enumerate}\end{prop}

Next for the sake of illustrating the above results, we give a
detailed analysis of $h_0$ in the radial case. Let us take
$\Omega=D(0,1)$ the unit disc in $\mathbb{R}^2$, $S_1=D(0,R)$ and
$S_2=D(0,1)\setminus \overline {D(0,R)}$ where $0<R<1$. In this case
$h_0$ is radially symmetric, $h_0(x)=h_0(|x|)$, and so it solves the
following ODE,
\begin{equation}\label{eq:h0(r)}
\left\{\begin{array}{l}
-h''(r)-r^{-1}h'(r)+h_0(r)=0\quad{\rm if}~0<r<R\,,\\
-h''(r)-r^{-1}h'(r)+a h_0(r)=0\quad{\rm if
}~R<r<1\,,\\
h_0(R_-)=h_0(R_+)\,,\quad h_0'(R_-)=\frac1 a h_0'(R_+)\,,\\
h_0'(0)=0\quad h_0(1)=1\,.
\end{array}\right.
\end{equation}
We look for a power series solution in the form,
$$h_0(r)=\sum_{n=0}^\infty a_n r^n\quad{\rm if}~0<r<R\,,\quad
h_0(r)=\sum_{n=0}^nb_n (r-R)^n\quad{\rm if}~R<r<1\,.$$ The sequences
$(a_n)$ and $(b_n)$ depend on $a$ and $R$ and all the terms are
expressed as functions of the term $a_0$, the two constants
$$\alpha=1+\sum_{k=1}^\infty\frac{R^{2k}}{2^k((k+1)!)^2}\,,\quad\beta=\sum_{k=0}^\infty(2k+2)
\frac{R^{2k+1}}{2^k((k+1)!)^2}\,,$$ and the sequence $(\gamma_n(a))$
defined recursively by,
$$\gamma_{-3}(a)=1\,,\quad\gamma_{-2}(a)=-\frac{1}R\,,\quad\gamma_{n+1}(a)=
-\frac1R\gamma_n(a)+a\gamma_{n-1}(a)\,,\quad\forall~n\geq-1\,.$$
Actually, we may verify that,
$$a_{2n+1}=0\,,\quad
a_{2n+2}=\frac{a_0}{2^n((n+1)!)^2}\quad\forall~n\in\mathbb N\,,$$
$$b_0=a_0\alpha\,,\quad b_1=a\,a_0\beta\,,\quad
b_{n+2}=\frac1{(n+2)!}\left(\gamma_{n-2}(a)b_1+a\gamma_{n-3}(a)b_0\right)\quad\forall~n\geq2\,,$$
and the term $a_0$ is expressed explicitly by,
$$a_0=\left[\alpha+
a\left(\beta(1-R)+\sum_{n=0}^\infty\frac1{(n+2)!}\left(\alpha\gamma_{n-2}(a)+\beta\gamma_{n-3}(a)\right)
(1-R)^n\right)\right]^{-1} \,,$$ provided that the sum in the r.h.s.
is finite. For instance, this is the case when $\frac12<R<1$ and
$a\to0_+$. Actually, in the limit $a\to0_+$, it holds that,
$$a_0\to \alpha^{-1}\,,\quad \frac{1-b_0}{a}\to \alpha^{-1}\left(\beta(1-R)+\frac{\beta-\alpha}{R}
\sum_{n=0}^\infty\frac{(-1)^n}{(n+2)!}\left(\frac1R-1\right)^n\right)=:c_0\,.$$
For $R\in(\frac12,1)$ chosen conveniently (close to $\frac12$), it
is easy to verify that $c_0>1-\alpha^{-1}$. Therefore, there exists
$a_0\in(0,1)$ sufficiently small (depending on $R\in(\frac12,1)$)
such that, for all $a\in(0,a_0)$, we have,
$$\lambda_2(a)<\lambda_1(a)\,,$$
where $\lambda_1(a)$ and $\lambda_2(a)$ are the critical constants
introduced in (\ref{eq:l0}). Coming back to the problem of vortex
nucleation for the G-L energy (\ref{reducedfunctional*}), we get in
light of Theorem~\ref{thm2} and Proposition~\ref{prop:w-l}, that
for a wide range of applied magnetic fields,
$$H=\lambda|\ln\varepsilon|(1+o(1))\quad(\varepsilon\to0)\,,\quad
\lambda\in(\lambda_2(a),\lambda_1(a))\,,$$ vortices exist and are
pinned in $S_2$. This result is in accordance with that obtained by
the second author in \cite{Kach-v1, Kach-v}, where nucleation of
vortices near the critical magnetic field (in the case of the disc)
is studied in details.

\section*{Acknowledgments}
Part of this work was carried out while the second author visited
Lebanese University. He wishes to thank A. Mneimneh and R. Talhouk
for their hospitality.

\end{document}